\RequirePackage{fix-cm}
\documentclass[smallextended]{svjour3}       
\smartqed  
\usepackage{amssymb} 

\journalname{Acta Applicandae Mathematicae}

\newtheorem{thm}{Theorem}

\def \Eta{\rm H} 
\def \D {\hbox{d}}
\def \Elemsimp {\mathop{\rm H}\nolimits}
\def \mod#1{\vert #1 \vert}

\def \Re#1{{\rm Re} #1} 
\def \Im#1{{\rm Im} #1} 
\def \csi{\kappa_{\rm i}} 
 
\def \ex{e_1} 
\def \ey{e_0} 
\def \dlogA{B} 
\def\Arbzero{M_{0}}
\def\Arbone{M_{1}}

\def\re{\rho}
\def \jmax{J}

\begin{document}

\title{Meromorphic traveling wave solutions of 
the complex cubic-quintic Ginzburg-Landau equation\thanks{
Part of this work was supported by 
RGC grant HKU 703807P and by France-Hong Kong grant F-HK39/11T. 
RC gladfully acknowledges the support of MPIPKS Dresden.
}} 


\titlerunning{Meromorphic solutions of CGL5} 

\author{Robert Conte        \and
        Tuen-Wai Ng 
}


\institute{Robert Conte \at
1. LRC MESO, 
Centre de math\'ematiques et de leurs applications (UMR 8536) 
\\ et CEA-DAM, \'Ecole normale sup\'erieure de Cachan, 
\\ 61, avenue du Pr\'esident Wilson, F--94235 Cachan Cedex, France.
\smallskip
\\ \noindent 2.
Department of Mathematics,
The University of Hong Kong,
Pokfulam Road.
       \email{Robert.Conte@cea.fr} 
           \and
           Tuen-Wai Ng \at
Department of Mathematics,
The University of Hong Kong,
Pokfulam Road.    
       \email{ntw@maths.hku.hk}              
}

\date{Received: 6 January 2012 / Revised: 27 March 2012}

\maketitle

\begin{abstract}
We look for singlevalued solutions of the squared modulus $M$
of the traveling wave reduction
of the complex cubic-quintic Ginzburg-Landau equation.
Using Clunie's lemma,
we first prove that any meromorphic solution $M$ 
is necessarily elliptic or degenerate elliptic.
We then give the two canonical decompositions of the new elliptic solution
recently obtained by the subequation method.
\keywords{Elliptic solutions \and
complex quintic Ginzburg-Landau equation}


\PACS{PACS 02.30.H  \and PACS 02.30.-f}


\subclass{MSC 30D30 \and MSC 33E05 \and 34A05}
\end{abstract}

\section{Introduction. The CGL5 and CGL3 equations}
\label{sectionIntro}

When a system is governed by an autonomous nonlinear algebraic 
partial differential equation (PDE),
it frequently admits permanent profile structures
such as fronts, pulses, sinks, etc
\cite{vS2003},
and usually these profiles are mathematically some 
singlevalued solution
of the traveling wave reduction $(x,t) \to x-ct$
of the PDE to an ordinary differential equation (ODE).

When the field is a slowly varying complex amplitude $A$,
the simplest equation involving time evolution, dispersion, nonlinearity and forcing
is the one-dimensional complex Ginzburg-Landau equation
\begin{eqnarray}
& &  {\hskip -18.0 truemm}
i A_t +p A_{xx} +q \mod{A}^2 A +r \mod{A}^4 A -i \gamma A =0,\
(A,p,q,r) \in \mathbb{C},\ \gamma \in \mathbb{R}. 
\label{eqCGL5}
\end{eqnarray}

We only consider in this class the equations 
which have the worst singularity structure,
the cubic one (CGL3, $r=0$, $\Im(q/p)\not=0$)
and the cubic-quintic one (CGL5, $\Im(r/p)\not=0$).
For a summary of results, 
see the reviews \cite{AK2002,vS2003}. 

CGL5 depends on seven real parameters.
Its travelling wave reduction
\begin{eqnarray}
& & 
A(x,t)=\sqrt{M(\xi)} e^{i(\displaystyle{-\omega t + \varphi(\xi)})},\
\xi=x-ct,\
(c,\omega,M,\varphi) \in {\mathbb R},
\label{eqCGL35red}
\\
& &
{\hskip -10.0 truemm}
 \frac{M''}{2 M} -\frac{{M'}^2}{4 M^2} + i \varphi''- {\varphi'}^2
        + i \varphi' \frac{M'}{M}
- i \frac{c}{2p} \frac{M'}{M}
+ \frac{c}{p} \varphi' + \frac{q}{p} M
+ \frac{r}{p} M^2
\nonumber \\ & &
+ \frac{\omega - i \gamma}{p}
=0,
\label{eqCGL5ReducComplex}
\end{eqnarray}
depends on eight real parameters, denoted
$e_r,e_i,d_r,d_i,s_r,s_i,g_r,g_i$,
\begin{eqnarray}
& &
e_r + i e_i = \frac{r}{p},\
d_r + i d_i = \frac{q}{p},\
s_r - i s_i = \frac{1}{p},\
\nonumber \\ & &
g_r + i g_i=\frac{\gamma + i \omega}{p}
 + \frac{c^2 s_r}{4} (2 s_i + i s_r).
\end{eqnarray}
The fact of taking account of the phase invariance \cite{vS2003}
\begin{eqnarray}
& &
\varphi' = \psi + \Re \frac{c}{2 p},
\label{eqdefpsi}
\end{eqnarray}
further reduces them to 
the seven real parameters $e_r,e_i,d_r,d_i,g_r,g_i,c s_i \equiv \csi$.

This third order system (\ref{eqCGL5ReducComplex}) can be written 
either 
as a real two-component rational system in the real variables $(M,\psi)$,
\begin{eqnarray}
& &
\left\lbrace
\begin{array}{ll}
\displaystyle{
\frac{M''}{2 M} -\frac{{M'}^2}{4 M^2} - \csi \frac{M'}{2 M} 
- \psi^2+ e_r M^2+ d_r M + g_i =0,
}\\ \displaystyle{
\psi' + \psi \frac{M'}{M} - \csi \psi + e_i M^2 + d_i M - g_r =0,
}
\end{array}
\right.
\label{eqCGL5ReducRealSystem}
\end{eqnarray}
\hfill\break\noindent
or, by elimination of $\psi$, 
as a real third order second degree ODE in $M$
\cite{Klyachkin} 
\begin{eqnarray}
& &
\psi = \frac{2 \csi G - G'}{2 M^2 (e_i M^2 + d_i M - g_r)},\ 
\psi^2=\frac{G}{M^2},
\label{eqCGL5Phiprime}
\\
& &
(G'-2 \csi G)^2 - 4 G M^2  (e_i M^2 + d_i M - g_r)^2=0,\ 
\label{eqCGL5Order3}
\\
& &
G=\frac{1}{2} M M'' - \frac{1}{4} M'^2
  -\frac{\csi}{2} M M' + e_r M^4 + d_r M^3 + g_i M^2.
\end{eqnarray}
 
The purpose of this work is to show that,
for all values of the seven parameters,
all meromorphic particular solutions $M$ 
of CGL3 and CGL5 belong to class $W$ (like Weierstrass),
defined as consisting of 
elliptic functions and their successive degeneracies,
i.e.: 
elliptic functions (meromorphic doubly periodic), 
rational functions of one exponential $e^{k \xi}$, $k\in\mathbb{C}$,
rational functions of $\xi$.
The assumption $M$ meromorphic 
implies the same property for
the variables $M'/M,\psi$ and the logarithmic derivative
of the complex amplitude $A e^{i \omega t}$,
\begin{eqnarray}
& & 
\dlogA:=\frac{\D}{\D \xi} \log \left(A e^{\displaystyle{i \omega t - i c s_r \xi/2}}\right)
 =\frac{M'}{2 M} + i \psi.
\label{eqdlogA}
\end{eqnarray}

In case $q \csi=0$,
the real system (\ref{eqCGL5ReducRealSystem})
displays a parity invariance,
\begin{eqnarray}
& q=0:\ & (M,\psi,\xi) \to (-M,\psi,\xi),
\label{eqInvariance_q}
\\
& \csi=0:\ & (M,\psi,\xi) \to (M,-\psi,-\xi).
\label{eqInvariance_csi}
\end{eqnarray}
This paper presents results complementary to those of \cite{ConteNgCGL5traps}.
Section \ref{sectionPZ} recalls the 
singularities
of $M$ and $\psi$. 
In section \ref{sectionClunie},
we prove that,
for all values of the seven real parameters,
any meromorphic solution $M$ of (\ref{eqCGL5ReducRealSystem})
is in class $W$.
In section \ref{sectionSubeqMethod},
we recall a method to obtain 
all the elliptic or degenerate elliptic solutions $M$,
and present the first order ODE for $M(\xi)$ 
characterizing the recently obtained elliptic solution \cite{ConteNgCGL5traps}.
Finally, section \ref{AppendixLaurent4} is devoted to
the construction of canonical expressions
to represent this elliptic solution.

\section{Movable singularities of CGL3 and CGL5} 
\label{sectionPZ}

Our interest is to count the number of distinct Laurent series for $M$ and $\psi$.
The results, obtained in \cite{CT1989} for CGL3 
and in \cite{MCC1994,ConteNgCGL5traps} for CGL5,
are the following.

A first set of poles $\chi_1=\xi-\xi_1 \to 0$ for both $M$ and $\psi$
arises by balancing $A_{xx}$ and the highest nonlinearity
($\mod{A}^2 A$ for CGL3, $\mod{A}^4 A$ for CGL5),
\begin{eqnarray}
& & {\hskip -11.0 truemm}
\hbox{CGL3}:\
A            \sim A_0 \chi_1^{-1+i \alpha},\ 
\overline{A} \sim A_0 \chi_1^{-1-i \alpha},\ 
(1-i \alpha) (2-i \alpha) p + A_0^2 q=0,
\label{eqleadingA3}
\\ & & {\hskip -11.0 truemm}
\hbox{CGL5}:\
A            \sim A_0 \chi_1^{-\frac{1}{2}+i \alpha},\ 
\overline{A} \sim A_0 \chi_1^{-\frac{1}{2}-i \alpha},\ 
(\frac{1}{2}-i \alpha) (\frac{3}{2}-i \alpha) p + A_0^4 r=0,
\label{eqleadingA5}
\end{eqnarray}
they define two real values $\alpha$
and two real values $A_0^2$ (CGL3) and $A_0^4$ (CGL5),
\begin{eqnarray}
& & {\hskip -12.0 truemm}
\hbox{CGL3}:\ 
d_i \alpha^2 -3 d_r -2 d_i=0,\ A_0^2=3 \alpha/d_i, 
\\ & & {\hskip -12.0 truemm}
\hbox{CGL5}:\ 
4 e_i \alpha^2 -8 e_r -3 e_i=0,\  A_0^4=2 \alpha/e_i,
\end{eqnarray}
and the Fuchs indices are $-1,0$ and two irrational numbers
\begin{eqnarray}
& & {\hskip -12.0 truemm}
\hbox{CGL3}:\ \hbox{indices =} 
-1,0, (7 \pm \sqrt{1-24 \alpha^2})/2, 
\\ & & {\hskip -12.0 truemm}
\hbox{CGL5}:\ \hbox{indices =} 
-1,0, (5 \pm \sqrt{1-32 \alpha^2})/2. 
\end{eqnarray}

For CGL3, $M$ presents double poles and $\psi$ simple poles,
\begin{eqnarray}
& & {\hskip -12.0 truemm}
M=m_0 \chi_1^{-2}
\left[
1+\frac{\csi}{3} \chi_1
+\mathcal{O}(\chi_1^2)
\right],\ m_0=A_0^2,
\label{eqCGL3Laurent_poles_M}
\label{eqCGL3_poles_M} 
\\ & & {\hskip -12.0 truemm}
\psi=\frac{d_i m_0}{3} \chi_1^{-1}
\left[
1+\frac{\csi}{6} \chi_1
+\mathcal{O}(\chi_1^2)
\right],\ 
\label{eqCGL3Laurent_poles_psi}
\end{eqnarray}
and for CGL5 $M$ and $\psi$ present simple poles,
\begin{eqnarray}
& & {\hskip -12.0 truemm}
M=m_0 \chi_1^{-1}
\left[
1+\left(\frac{\csi}{4}+\frac{2 d_r m_0-2 e_i d_i m_0^3}{4(1+e_i^2 m_0^4)}\right) \chi_1
+\mathcal{O}(\chi_1^2)
\right],\ m_0=A_0^2,
\label{eqCGL5Laurent_poles_M}
\label{eqCGL5_poles_M} 
\\ & & {\hskip -12.0 truemm}
\psi=\frac{e_i m_0^2}{2} \chi_1^{-1}
+ \frac{e_i m_0^2}{8}\csi + m_0 \frac{4 d_i +5 e_i d_r m_0^2- e_i^2 d_i m_0^4}{4(1+e_i^2 m_0^4)}
+\mathcal{O}(\chi_1),
\label{eqCGL5Laurent_poles_psi}
\end{eqnarray}
in which both invariances (\ref{eqInvariance_q})--(\ref{eqInvariance_csi})
require changing $m_0$ to $-m_0$.
The number of distinct Laurent series $M$ near $\chi_1$ 
is two (CGL3) or four (CGL5),
and the number of series $\psi$ is two (CGL3), four (CGL5 $q\not=0$) or two (CGL5 $q=0$).

A second set of singularities is
easier to compute from the system (\ref{eqCGL5ReducRealSystem}),
\begin{eqnarray}
& & {\hskip -3.0truemm}
\frac{1}{M}=\frac{1}{\Arbzero}\chi_2^{-1}
\left[
1 + \Arbone \chi_2
+\left\lbrace \Arbone^2 + \csi \Arbone - \frac{j}{3} g_r + \frac{2}{3} g_i  
\right\rbrace \chi_2^2
+\mathcal{O}(\chi_2^3)
\right],
\label{eqCGL5Laurent_zeroes_M}
\\ & & {\hskip -3.0truemm}
\psi= \frac{j}{2} \chi_2^{-1}
\left[
1 
+ \left(\csi + \Arbone\right) \chi_2
+\left\lbrace \Arbone^2 + 2 \csi \Arbone + \frac{2}{3} g_i - \frac{4 j}{3} g_r +\frac{5}{6}\csi^2
\right\rbrace \chi_2^2
\right.
\nonumber \\ & & {\hskip -3.0truemm}
\left. \phantom{12345}
+\frac{1}{2}\left\lbrace 
 (g_r + j g_i)\csi
 + \frac{3 j \csi^3}{4} 
-\frac{(3 d_i -j d_r)\Arbzero}{4} 
+ \frac{(11 j \csi^2 + 4 g_r + 4 j g_i)\Arbone}{4} 
\right. \right.
\nonumber \\ & & {\hskip -3.0truemm}
\left. \left. \phantom{123456789}
+ 3 j \csi \Arbone^2
+ j \Arbone^3
\right\rbrace \chi_2^3
+\mathcal{O}(\chi_2^4)\right],
\label{eqCGL5Laurent_poles_psi_at_zeroes_M}
\end{eqnarray}
in which $\Arbzero$, $\Arbone$ are arbitrary constants,
and $j^2=-1$. 
Invariances (\ref{eqInvariance_q})--(\ref{eqInvariance_csi})
require changing $\Arbzero$ to $-\Arbzero$,
with $\Arbone=0$ when $\csi=0$.
This defines 
either $2 N$ (when $q\not=0$) 
or     $  N$ (when $q=0$) simple poles of $\psi$, with $N$ an undetermined integer.
A direct study \cite{ConteNgCGL5traps} of 
the third order ODEs for $M$ and $\psi$ shows that
neither $M$ nor $\psi$ admit other movable poles.

For CGL3 (resp.~CGL5), $M$ and $\psi M$ admit two (resp.~four) Laurent series.

\section{Results from Clunie's lemma}
\label{sectionClunie}

For convenience, $\xi$ will be denoted as $z$ in this section only. 
We shall prove 
\begin{thm}
For all values of the 
constants $p,q,r,\gamma,c,\omega$, 
all meromorphic traveling wave solutions $M$ 
of CGL3 and CGL5 equations belong to class $W$. 
\label{thm1}
\end{thm}

The method we use is a refinement of Eremenko's method 
developed in \cite{Eremenko1982} as well as \cite{EremenkoKS,ELN,ConteNgODE3},
based on the local singularity analysis 
of the solutions of the given differential equation
and on the zero distribution and growth rate of their meromorphic 
solutions by using Nevanlinna theory.

Several partial results have been previously obtained \cite{MC2003,Hone2005,VernovCGL5}
for finding solutions of (\ref{eqCGL5Order3}),
but they are incomplete and Theorem \ref{thm1} settles the question.
These previous results are the following.
\begin{enumerate}

\item
For CGL3, when $d_r\not=0$, 
all solutions belonging to class $W$ have been found \cite{MC2003}:
there are six distinct solutions which are rational functions in one exponential function 
and 
there is no elliptic solution.

\item
For CGL3, when $\csi\not=0$, there exists no 
elliptic solution \cite{Hone2005}.

\item
For CGL5, when $\csi=0$, 
there exists exactly one 
elliptic solution \cite{VernovCGL5}.

\end{enumerate}

Let us recall a definition. For a differential polynomial of $f$,
$$P(z,f)=\sum a_j f^{j_0}(f')^{j_1}\cdots(f^{(k)})^{j_k},$$ 
where $j = (j_0, . . . , j_k)$ is a multi-index and $f$ and $a_j$ are meromorphic functions, 
the sum $j_0 +. . .+j_k$ is called the \textit{degree} of the monomial 
$a_jf^{j_0}(f')^{j_1}\cdots(f^{(k)})^{j_k}$. 
The \textit{total degree} of $P(z,f)$ is defined as
the maximum of the degrees of its monomials.

We shall assume the readers are familiar with the terminology 
and results of Nevanlinna theory 
\cite{Hayman1964,Laine-book,Nevanlinna-book} 
(see \cite{EremenkoKS} for a quick introduction). 
Here, we recall some basic notations of Nevanlinna theory. 
Let $f$ be a non-constant meromorphic function on the open disc $D(r)=\{z: |z|<r\}$ where $r$. 
Denote the number of poles of $f$ on the closed disc $\overline{D}(r)$ by $n(r,f)$, 
counting multiplicity. 
Define the \textit{integrated counting function} $N(r,f)$ by
$N(r,f)=n(0,f)\log r +\int^{r}_{0} \left[n(t,f)-n(0,f)\right] \frac{dt}{t}$ 
and the \textit{proximity function} $m(r,f)$ by
$m(r,f)=\int^{2\pi}_{0}\log^{+} f(re^{i\theta})\frac{d\theta}{2\pi}$, 
where $\log^{+}x=\max {\left\{0,\log x\right\}}$. 
Finally, the \textit{Nevanlinna characteristic function} $T(r,f)$ 
is defined by $T(r,f)=m(r,f)+N(r,f)$ and we let $S(r,f)$ be a term such that 
$\frac{S(r,f)}{T(r,f)} \rightarrow 0 \,,$ as $\, r\rightarrow +\infty$.

To prove Theorem \ref{thm1}, we make use of the well known Clunie's Lemma,

\begin{lemma}  \cite[2.4.2]{Laine-book} 
Let $f$ be a transcendental meromorphic solution of 
\begin{eqnarray}
& &
f^n P(z,f)=Q(z,f),
\label{eqClunie_input}
\end{eqnarray}
with $n$ a positive integer, $P$ and $Q$ differential polynomials of $f$ 
with meromorphic coefficients 
$a_\lambda$ such that  $m(r,a_\lambda)=S(r,f)$. 
If the total degree of $Q$ 
is less than or equal to $n$, then 
\begin{eqnarray}
& &
m(r,P(z,f))=S(r,f).
\label{eqClunie_output}
\end{eqnarray}

\label{lemma_Clunie}
\end{lemma} 

Actually, all we need is the following corollary of Clunie's Lemma.
\begin{corollary}  Let $f$ be a transcendental meromorphic solution of the ODE 
\begin{eqnarray}
& &
f^{n+1}=Q(z,f),
\label{DE_input}
\end{eqnarray}
with $n$ a positive integer and $Q$ a differential polynomial of $f$ 
with meromorphic coefficients 
$a_\lambda$ such that  $m(r,a_\lambda)=S(r,f)$. 
If the total degree of $Q$ is less than or equal to $n$, 
then $f$ must have infinitely many poles. 
\label{corollary_Clunie}
\end{corollary} 
\noindent
\proof of Corollary \ref{corollary_Clunie}.

Taking $P(z,f)=f$ in Lemma \ref{lemma_Clunie}, 
we conclude that $m(r,f)=S(r,f)$, and therefore $(1-o(1))T(r,f)=N(r,f)$. 
Assume that $f$ has finitely many poles. 
Then $N(r,f)= O(\log r)$, and therefore $T(r,f)=O(\log r)$, 
which is impossible since $f$ is transcendental. \qed

\noindent
\proof of Theorem \ref{thm1}

Let $M$ be a solution of (\ref{eqCGL5Order3}) which is meromorphic in the complex plane.
If $M$ is rational, then we are done.
So suppose $M$ is transcendental and let us prove that $M$ has infinitely many poles. 
We first rewrite the second equation of (\ref{eqCGL5ReducRealSystem}) as 
\begin{eqnarray}
\begin{array}{ll}
(\psi M)' - \csi (\psi M) + e_i M^3 + d_i M^2 - g_rM =0.
\end{array}
\label{eqCGL5ReducRealSubSystem}
\end{eqnarray}
It follows easily that if $\psi$ has infinitely many poles, then so does $M$. 

We first show that if $\psi$ is transcendental, then $\psi$ has infinitely many poles and hence so does $M$. One can build an ODE for $\psi(z)$ via the elimination of $M$
between the system (\ref{eqCGL5ReducRealSystem}). This third order ODE is given as follows,
\begin{eqnarray}
& &
\hbox{CGL5 }:\ e_i^2 (5 e_i^2+e_r^2)^2 \psi^{20}=Q(z,\psi),
\label{eqClunie_form_CGL5psi}
\\
& &
\hbox{CGL3 }:\ d_i (3 d_i^2 + d_r^2) \psi^{10}=Q(z,\psi), 
\label{eqClunie_form_CGL3psi}
\end{eqnarray}
where the differential polynomial $Q(z,\psi)$ has the total degree $19$ (CGL5 case) or $9$ (CGL3 case). 
\medskip
Applying Corollary \ref{corollary_Clunie} to
(\ref{eqClunie_form_CGL5psi}) and (\ref{eqClunie_form_CGL3psi}),
we conclude that $\psi$ and therefore $M$ must have infinitely many poles. 

Now suppose $\psi$ is rational, then it is well known that 
$T(r,\psi)=O(\log r)$ and $T(r,\psi')=O(\log r)$,
and therefore $m(r,\psi)=S(r,f)$ and so is $m(r,\psi')$ since $M$ is transcendental. 
Now (\ref{eqCGL5ReducRealSubSystem}) can be written as 
$e_i M^3= -\psi' M - \psi M + (\csi \psi) M - d_i M^2 + g_rM$
 and applying Corollary \ref{corollary_Clunie} to it, we again conclude that $M$ has infinitely many poles.

Secondly, knowing that the transcendental meromorphic solution $M$
has infinitely many poles, let us prove that it is a periodic function. 
By the local singularity analysis (section \ref{sectionPZ}),
if $z_0$ is a pole of $M$, 
CGL3 (resp.~CGL5) admits exactly two (resp.~four) Laurent series $M$ 
with poles at $z=z_0$ obeying the ODE (\ref{eqCGL5Order3}).
Now let $z_j, j=1,2,3,\cdots$ be the poles of $M(z)$;
the functions $w_j(z) = M(z+z_j-z_0)$ are then meromorphic 
solutions of the ODE (\ref{eqCGL5Order3}) with a pole at $z_0$,
therefore some of them must be equal. 
Consequently, $M$ is a periodic function. 

Without loss of generality, we assume that $M$ has a period of $2\pi i$. 
Let $D=\{z:0\le {\rm Im} z <2\pi\}$. 
If $M$ has more than three (CGL3) or five (CGL5) poles in $D$, then by the previous argument
we conclude that $M$ is periodic in $D$ 
and therefore is indeed an elliptic function and we are done.

Now suppose $M$ has at most two (CGL3) or four (CGL5) poles in $D$.
Since $M$ is a periodic function with period $2\pi i$, we have
$N(r,M)= O(r),$ as $r\to\infty$. 
It follows from  $(1-o(1))T(r,M)=N(r,M)$ that $T(r,M)=O(r)$. 
By Nevanlinna's First Fundamental Theorem, we know that 
for any $a\in\mathbb{C}$, $N(r,1/(M -a)) =O(r)$ as $r\to\infty$. 
By the periodicity of $M$, we conclude that $M$ takes each $a$ finitely many
times in $D$. 
Hence, the function $R(z)=M(\ln z)$ is a
single-valued analytic function in the punctured plane
$\{z:0<|z|<\infty\}$ and takes each $a$ finitely many times. 
It follows that $0$ is a removable singularity of $R$, and $R$ 
must then be a rational function. 
Therefore, $M(z)=R(e^z)$ belongs to class $W$. \qed

\section{A method to determine all solutions in class $W$} 
\label{sectionSubeqMethod}

Consider an $N$-th order autonomous algebraic ODE, 
\begin{eqnarray}
& &
E(u^{(N)},...,u',u)=0,\
'=\D / \D x, 
\label{eqODE}
\end{eqnarray}
admitting at least one Laurent series
\begin{eqnarray}
& &
u=\chi^p \sum_{j=0}^{+\infty} u_j \chi^j,\ \chi=x-x_0.
\label{eqLaurent}
\end{eqnarray}

There exists an algorithm \cite{MC2003} to find in closed form
all its elliptic or degenerate elliptic solutions.
Its successive steps are \cite{CM2009,CMBook}:
\begin{enumerate}
\item
Find the structure of movable singularities
(e.g., 4 families of simple poles).
For each subset of families (e.g.~2 families of simple poles) 
deduce the elliptic orders $m,n$ (e.g.~$m=2,n=4$) of $u,u'$
and perform the next steps.

\item
Compute slightly more than $(m+1)^2$ terms in the Laurent series.

\item
Define the first order $m$-th degree
subequation 
\begin{eqnarray}
& &
F(u,u') \equiv
 \sum_{k=0}^{m} \sum_{j=0}^{2m-2k} a_{j,k} u^j {u'}^k=0,\ a_{0,m}\not=0.
\end{eqnarray}
According to 
results of Briot-Bouquet and Painlev\'e \cite{CM2009},
any solution of (\ref{eqODE}) in class $W$
\textit{must} obey such an ODE,
called a ``subequation'' 
because (next step) it admits (\ref{eqODE}) as a differential consequence.

\item
Require each Laurent series (\ref{eqLaurent}) to obey $F(u,u')=0$,
\begin{eqnarray}
& & {\hskip -10.0 truemm}
F \equiv \chi^{m(p-1)} \left(\sum_{j=0}^{\jmax} F_j \chi^j
 + {\mathcal O}(\chi^{\jmax+1})
\right),\
\forall j\ : \ F_j=0.
\label{eqLinearSystemFj}
\end{eqnarray}
and solve this \textbf{linear overdetermined} system for $a_{j,k}$.

\item
Integrate each resulting ODE $F(u,u')=0$.
\end{enumerate}

A similar method has later been developed \cite{DK2011first},
which also takes advantage of the Laurent series
and directly searches for a canonical closed form representation
of the elliptic solutions and their degeneracies.

Theorem \ref{thm1} 
implies that the subequation method 
is indeed able to find all the meromorphic traveling wave solutions
of the CGL3 and CGL5 equations. 

For CGL5, the subequation method
has produced a new elliptic solution \cite{ConteNgCGL5traps},
characterized by the first order, fourth degree, genus one ODE
\begin{eqnarray}
& & {\hskip -10.0 truemm}
F_4 \equiv
{M'}^4 
-2 \csi M {M'}^3
 +\frac{72}{e_i} \ex {M'}^2 (e_i M^2 - 12 \ey)
 +\frac{2^4 3^8 \ex^4}{e_i^2}
\nonumber \\ & &  {\hskip -10.0 truemm}
 +\frac{648 \ex^2}{e_i^2} \left(288 \ey^2 +24 e_i \ey M^2 -e_i^2 M^4\right)
 - \frac{1}{3^4 e_i} M^2 \left(e_i M^2 -48 \ey\right)^3
=0,
\label{eqsub54Mexey}
\\ 
& & {\hskip -10.0 truemm}
\csi^2=48 \ex,\ g_r=36 \ey,\ e_r=d_r=d_i=0,\ g_i=-\frac{3}{16} \csi^2. 
\end{eqnarray}


\section{Integration of subequation (\ref{eqsub54Mexey})}
\label{AppendixLaurent4}

Let us first recall the differential equations of Weierstrass
\begin{eqnarray}
& &
{\wp'}^2= 4 \wp^3 - g_2 \wp - g_3,\
\zeta'=-\wp,\
(\log \sigma)'=\zeta.
\end{eqnarray}

Apart the representation as a rational function of $\wp$ and $\wp'$, 
\begin{eqnarray}
& & 
\frac{\hbox{polynomial}(\wp) + \hbox{polynomial}(\wp) \wp'}{\hbox{polynomial}(\wp)},
\label{eqsimple0}
\end{eqnarray}
elliptic functions have two main decompositions,
either as a sum 
\begin{eqnarray} 
& & 
C+\sum_{j=1}^{N} \left(
 r_j\ \zeta(\xi-a_j)+\sum_{k=0}^{M} c_{j,k} \wp^{(k)}(\xi-a_j)
\right),\
 \sum_{j=1}^{N} r_j=0,\   
\label{eqsimplezeta}
\end{eqnarray}
in which $C, r_j, a_j, c_{j,k}$ are complex constants ($a_j$ distinct),
or as a quotient of two products of an equal number of entire  functions $\sigma$,
\begin{eqnarray}
& & 
\hbox{constant } \prod_{j=1}^{P} \frac{\sigma(\xi-\alpha_j)}{\sigma(\xi-\beta_j)},\ 
\sum_{j=1}^{P} \alpha_j-\beta_j=0,\  
\label{eqsimplesigma}
\end{eqnarray}
in which $\alpha_j, \beta_j$ are not necessarily distinct complex constants.

To obtain the complex amplitude $A$, which is not elliptic,
one can either compute the couple $(M,\psi)$ then perform the quadrature $\int \psi\ \D \xi$,
or compute the logarithmic derivative $\dlogA$, Eq.~(\ref{eqdlogA}),
then perform the quadrature $\int \dlogA\ \D \xi$.

By elimination with (\ref{eqCGL5Phiprime}),
one first deduces the real subequation for $\psi$,
\begin{eqnarray}
& & 
  \csi {\psi'}^4 
- 4 \csi {\psi'}^3 \left(\csi \psi + 24 \ey\right)
\nonumber \\ & & 
+ 8 {\psi'}^2 \Big(-\csi (27 \ex^2-324 \ey^2)
              +1440 \ex \ey \psi
              +27 \csi \ex \psi^2
              +16 \ey \psi^3
              +\frac{1}{3} \csi \psi^4\Big)             
\nonumber \\ & & 
 +16 \Big( 
              -\frac{1}{3} \csi \psi^8
              -\frac{32}{3} \ey \psi^7
              -26 \csi \ex      \psi^6
              -1632 \ex \ey \psi^5                 
             - \left(477 \ex^2+552 \ey^2 \right)\psi^4
     \Big.
\nonumber \\ & & 
             - 288 \left(165 \ex^2+4 \ey^2\right)\psi^3                 
             +\csi \left(2106 \ex^2-31320 \ey^2\right)\psi^2
\nonumber \\ & & 
\Big. 
             +2^7 3^6 \left(\ex^2-4 \ey^2\right) \ex \ey \psi    
             +243    \left(-9 \ex^4 +56 \ex^2 \ey^2 -144 \ey^4\right)  
\Big)=0,
\label{eqsub54psi}
\end{eqnarray}
then the complex subequation 
for 
$\dlogA$ as defined Eq.~(\ref{eqdlogA}),
\begin{eqnarray}
& & 
\left(2 \dlogA' + \csi \dlogA + 24 i \ey \right) 
\left(  \dlogA' - \csi \dlogA - 24 i \ey \right)^2
\nonumber \\ & & 
+ 2^{-11} 
\left(16 (4 \dlogA^3-3 \csi \dlogA^2) -9 (\csi^2+64 i \ey) (4 \dlogA + \csi) \right)^2
=0.   
\label{eqsub54dlogA}
\end{eqnarray}

The degree of subequation (\ref{eqsub54psi}) drops from four to two
when $\csi=0$.
As to (\ref{eqsub54dlogA}),
it has degree three
and therefore belongs to the so-called trinomial type integrated by
Briot and Bouquet \cite[\S 250--251 p.~395]{BriotBouquet}.

Let us derive decompositions (\ref{eqsimplezeta}) or (\ref{eqsimplesigma})  
for the solution of genus one equations  
(\ref{eqsub54Mexey}) for $M$,
(\ref{eqsub54psi}) for $\psi$ or 
(\ref{eqsub54dlogA}) for $\dlogA$.
Three steps are required. 
\begin{description}
\item --
The \textit{first step} is to represent the solution
as a rational function of 
$\wp(\xi-\xi_0)$ and $\wp'(\xi-\xi_0)$,
in which $\xi_0$ is arbitrary,
and to write it in the canonical form (\ref{eqsimple0}).
Because of the existence of an addition formula for $\wp$,
\begin{eqnarray}
& &
\forall x_1,x_2:\
\wp(x_1+x_2)+\wp(x_1)+\wp(x_2)=\frac{1}{4}
\left(\frac{\wp'(x_1)-\wp'(x_2)}{\wp(x_1)-\wp(x_2)}\right)^2,
\label{eqWeierstrassAddition}
\end{eqnarray}
such a canonical form (\ref{eqsimple0}) is not unique,
and general algorithms may yield
messy expressions by performing a shift on $\xi_0$.
For instance, 
with the (otherwise powerful) command
\verb+Weierstrassform+ \cite{MapleAlgcurves} of the computer algebra language Maple \cite{Maple},
which applies to any genus one equation,
the ODE
\begin{eqnarray}
& &
{u'}^2=u^4 - u^3 + u^2 + u +7
\end{eqnarray}
is integrated as the second degree rational function
\begin{eqnarray}
& &
u=3 \frac{12 \wp + 24 \sqrt{7} \wp'}{144 \wp^2 - 24 \wp -251},\ g_2=\frac{22}{3},\ g_3=\frac{277}{432},\ 
\end{eqnarray}
while a first degree rational function $c_0+c_1/(\wp - c_2)$ is sufficient.
The same occurs with the algorithm of Briot and Bouquet 
\cite[\S 250--251 p.~395]{BriotBouquet} to integrate
binomial or trinomial equations:
with (\ref{eqsub54dlogA}),
instead of yielding a second degree rational function  
rational in the fixed constants
(see (\ref{eqCGL5dlogAwp}) below),
it yields a third degree rational function algebraic in the fixed constants.
Consequently, the practical method used here is
to determine the smallest degrees of the three polynomials in (\ref{eqsimple0}),
then their coefficients by identification.
One thus finds for the solution $M$ of (\ref{eqsub54Mexey}),
\begin{eqnarray}
& & {\hskip -5.0 truemm}
\left\lbrace
\begin{array}{ll}
\displaystyle{
M=\frac{8 N_0 (3 \ex + 4 j \ey)(\wp-\ex) [3 \ex \wp^2 + 4(3 \ex^2 + 4 \ey^2) \wp + 4 \ex (3 \ex^2 + 5 \ey^2)]}
       {  24(3 \ex + 4 j \ey)(\wp-\ex) (\wp^2 - 2 \ex \wp - 8 \ex^2 - 12 \ey^2)
           + \csi P_2^{M} \wp'},\
}\\ \displaystyle{
P_2^{M}= \wp^2 + 4(\ex + j \ey) \wp + 4 \ex^2 - 4 j \ex \ey +12 \ey^2,\
j^2=-1,\
}\\ \displaystyle{
N_0^2=-\frac{324 j \ex^2}{e_i (3 \ex + 4 j \ey)},\
}\\ \displaystyle{
\wp:=\wp(\xi-\xi_0^{M},g_2,g_3),
}\\ \displaystyle{
{\wp'}^2=4 (\wp-\ex) (\wp^2 + \ex \wp + 7 \ex^2 + 12 \ey^2),\ 
}\\ \displaystyle{
g_2=-24 (  \ex^2 +  2 \ey^2),\
g_3=  4 (7 \ex^2 + 12 \ey^2) \ex.
}
\end{array}
\right.
\label{eqCGL5Mwp}
\end{eqnarray}
This expression will simplify greatly as (\ref{eqCGL5Msumzeta_g}). 
Because of the correspondence (\ref{eqCGL5Phiprime}),
the solution $\psi$ of (\ref{eqsub54psi}) 
involves the same square root $j$ of $-1$ as in (\ref{eqCGL5Mwp}).
When $\csi$ is nonzero this is
\begin{eqnarray}
& & {\hskip -18.0 truemm}
\left\lbrace
\begin{array}{ll}
\displaystyle{
\psi=-\frac{j \csi (9 \ex-4 j \ey)}{24 \ex}
 + \frac{P_2^{\psi} + Q_2^{\psi} \wp'}{12 \ex (3 \ex \wp+15 \ex^2+16 \ey^2) ((\wp + 2 \ex)^2+3 (3 \ex+4 j \ey)^2 )},\
}\\ \displaystyle{
P_2^{\psi}= -j \csi (3 \ex+4 j \ey) [(3 \ex+2 j \ey) ((9 \ex-4 j \ey) \wp^2+2 (-9 \ex-44 j \ey) \ex \wp)
}\\ \displaystyle{\phantom{1234567} 
     -945 \ex^4-1434 j \ex^3 \ey-1192 \ey^2 \ex^2-1440 j \ey^3 \ex-384 \ey^4],\ 
}\\ \displaystyle{
Q_2^{\psi}= 9 j \ex (\ex (\wp^2+22 \wp \ex+24 j \wp \ey)+121 \ex^3+48 \ex \ey^2+192 j \ex^2 \ey+128 j \ey^3), 
}\\ \displaystyle{
\wp:=\wp(\xi-\xi_0^{M},G_2,G_3),
}\\ \displaystyle{
{\wp'}^2=4 (\wp+2 \ex) (\wp^2-2 \ex \wp-35 \ex^2-48 \ey^2),
}\\ \displaystyle{
G_2=12 (13 \ex^2+16 \ey^2),\
G_3= 8 (35 \ex^2+48 \ey^2) \ex,
}
\end{array}
\right.
\label{eqCGL5psiwp}
\end{eqnarray}
while for $\csi=0$ it is
\begin{eqnarray}
& & {\hskip -18.0 truemm}
\csi=0:\
\psi=\sqrt{6 j \sqrt{3} \ey} \left(1+\frac{8 j \sqrt{3} \ey}{\wp(\xi-\xi_0^{M},192 \ey^2,0) - 4 j \sqrt{3} \ey}\right),\ 
\label{eqCGL5psiwp0}
\end{eqnarray}
or simply (but with yet another $g_2$),
\begin{eqnarray}
& & {\hskip -18.0 truemm}
\csi=0:\
\psi=\frac{\sqrt{3}}{2} \sqrt{\wp(\xi-\xi_0^{M},-768 \ey^2,0)}.
\label{eqCGL5psiwp0other}
\end{eqnarray}

Finally, the solution $\dlogA$ of (\ref{eqsub54dlogA}) is expressed as,
\begin{eqnarray}
& & {\hskip -18.0 truemm}
\left\lbrace
\begin{array}{ll}
\displaystyle{
\dlogA=\frac{\csi}{2} 
 - \frac{6 \csi (3 \ex+4 i \ey)^2+ (\wp+2 \ex + 3(3\ex+4 i \ey)) \wp'}{2 ((\wp+2 \ex)^2+3 (3 \ex+4 i \ey)^2)},
}\\ \displaystyle{
\wp:=\wp(\xi-\xi_0^{\dlogA},G_2,G_3),
}
\end{array}
\right.
\label{eqCGL5dlogAwp}
\end{eqnarray}
The properties of the above three expressions are:
all coefficients (except the global factor $N_0$) are rational in $(\csi,g_r)$,
the two different Weierstrass functions $\wp(.,g_2,g_3)$ and $\wp(.,G_2,G_3)$
are linked by a Landen transformation (Appendix A),
the relation (\ref{eqdlogA}) between $\dlogA,M,\psi$ holds true
when 
the square root $j$ of $-1$ is equal to $+i$
and the constant origins $\xi_0^{M},\xi_0^{\dlogA}$
are equal.

The degeneracy $\Delta=0$ implies $g_2=g_3=0$,
i.e.~it 
directly defines the reducible subequation 
$\left(3 M'\right)^4 - e_i^2 M^8=0$,
whose solutions are rational.
Because of this, 
even a four-family extension of the method used in 
\cite{MCC1994} would fail for CGL5.

\item --
The \textit{second step} is
to compute the partial fraction decomposition 
of the rational functions (\ref{eqsimple0}) of the variable $\wp$,
considering for a moment $\wp'$ as a parameter.
The rational function (\ref{eqCGL5Mwp}),
once converted to the canonical form (\ref{eqsimple0}),
admits four poles for each choice of $j$, 
and we characterize their affixes $\xi_{j,k}^{M}$  by choosing
the signs of $\wp'(\xi_{j,k}^{M},g_2,g_3)$ as follows,
\begin{eqnarray}
& & {\hskip -18.0 truemm}
\left\lbrace
\begin{array}{ll}
\displaystyle{
\wp(\xi_{j,k}^{M},g_2,g_3)=
\left(-3      +3 (j^k+\sqrt{3} j^{1-k}) \re+(-1)^k \re^2\right) \ex/6,\ 
}\\ \displaystyle{
\wp'(\xi_{j,k}^{M},g_2,g_3)=
( 9 j^k+3 ((-1)^{1+k}-j \sqrt{3}) \re+j^{2-k} \re^2) \ex\csi\re/36,\ 
}\\ \displaystyle{
\ex \re^2=3 j \sqrt{3} (3 \ex + 4 j \ey),\ j=\pm i,\ k=1,2,3,4.
}
\end{array}
\right.
\label{eqpolesM}
\end{eqnarray}

The rational function (\ref{eqCGL5psiwp}) 
admits one real pole and, for each $j=\pm i$,
two complex poles similarly characterized
as follows,
\begin{eqnarray}
& & {\hskip -8.0 truemm}
\csi\not=0:\
\left\lbrace
\begin{array}{ll}
\displaystyle{
\wp(\xi_{j}^{\psi},G_2,G_3)= - 5 \ex - 16 \ey^2 /(3 \ex),\ 
}\\ \displaystyle{
\wp'(\xi_{j}^{\psi},G_2,G_3)=
- 2 j \csi \ey (9 \ex^2 +16 \ey^2/(9 \ex^2),\ j=\pm i,
}
\end{array}
\right.
\label{eqpolepsi}
\\
& & {\hskip -8.0 truemm}
\left\lbrace
\begin{array}{ll}
\displaystyle{
\wp(\xi_{j,k}^{\psi},G_2,G_3)=
- 2 \ex + (-1)^k j \sqrt{3} (3 \ex + 4 j \ey),\ j=\pm i,\ k=0,1, 
}\\ \displaystyle{
\wp'(\xi_{j,k}^{\psi},G_2,G_3)=
(3 - (-1)^k j \sqrt{3})\csi (3 \ex + 4 j \ey)/2.
}
\end{array}
\right.
\label{eqpolespsi}
\end{eqnarray}
Finally, the two poles of (\ref{eqCGL5dlogAwp}) 
are just $\wp(\xi_{i,k}^{\psi},G_2,G_3)$.

\textit{Modulo} the periods of $\wp(.,G_2,G_3)$,
the affixes of these poles obey 
\begin{eqnarray}
& \xi_{ j,0}^{\psi}+\xi_{j,1}^{\psi}- \xi_j^{\psi}:
& \wp=-2 \ex,\ \wp'=0,\ \zeta=\Eta_1\  \hbox{ (half-period)},
\label{eqEta1}
\\
& \xi_{ j,0}^{\psi}+\xi_{j,1}^{\psi}:
& \wp=7 \ex,\ \wp'=- 6 j \csi \ey,\
\\
& \xi_{ j,0}^{\psi}-\xi_{j,1}^{\psi}:
& \wp=-5 \ex,\ \wp'= - 2 \sqrt{3} \csi \ey,\ 
\\
& \xi_{ j,k}^{\psi}+\xi_{-j,1-k}^{\psi}- \xi_j^{\psi}:
& \wp=\infty,\ \wp'=\infty\ \hbox{ (period)}.
\label{eq_relations_among_poles}
\end{eqnarray}
The Landen transformation maps
$\wp(\xi_{j,2}^{M},g_2,g_3)$ and $\wp(\xi_{j,4}^{M},g_2,g_3)$ 
to \hfill\break\noindent
$\wp(\xi_{j,0}^{\psi},G_2,G_3)$,
and maps
$\wp(\xi_{j,1}^{M},g_2,g_3)$ and $\wp(\xi_{j,3}^{M},g_2,g_3)$ 
to $\wp(\xi_{j,1}^{\psi},G_2,G_3)$.

Expressions (\ref{eqCGL5Mwp}), (\ref{eqCGL5psiwp}), (\ref{eqCGL5dlogAwp})
thus evaluate to the sum 
\begin{eqnarray}
& & 
\hbox{constant } +\sum_{j=1}^{4} 
\frac{\hbox{constant } + \hbox{constant } \wp'(\xi-\xi_0)}{\wp(\xi-\xi_0) - \wp(\xi_j)}.
\label{eqsimpl1}
\end{eqnarray}

\item --
The \textit{third step} is, 
using the classical identities 
\begin{eqnarray}
& & {\hskip -18.0 truemm}
\forall u,v:\
\left\lbrace
\begin{array}{ll}
\displaystyle{
\zeta(u+v)+\zeta(u-v)-2 \zeta(u)=\frac{ \wp'(u)}{\wp(u)-\wp(v)},
}\\ \displaystyle{
\zeta(u+v)-\zeta(u-v)-2 \zeta(v)=\frac{-\wp'(v)}{\wp(u)-\wp(v)},
}
\end{array}
\right.
\label{eqzetaformulae}
\end{eqnarray}
to convert (\ref{eqsimpl1}) into a finite sum of $\zeta$ functions.
\end{description}

The result for (\ref{eqCGL5Mwp}),
\begin{eqnarray}
& & {\hskip -6.0 truemm}
\forall \csi:\
M=\frac{3^{1/4}}{\sqrt{-e_i}} \sum_{k=1}^4
  j^{k-1} \left(\zeta(\xi-\xi_{j,k}^{M},g_2,g_3) +\zeta(\xi_{j,k}^{M},g_2,g_3)\right),\ j^2=-1,
\label{eqCGL5Msumzeta_g}
\end{eqnarray}
clearly displays the four simple poles.

As to the three simple pole variables $\dlogA, \psi, M'/M$,
their decompositions evaluate to
(we abbreviate $\zeta(.,G_2,G_3)$ to $\zeta(.)$ )
\begin{eqnarray}
& & {\hskip -18.0 truemm}
\left\lbrace
\begin{array}{ll}
\displaystyle{
\forall \csi:\
\frac{\D}{\D \xi} \log \left(A e^{\displaystyle{i \omega t - i \frac{c s_r}{2} \xi}}\right)=\frac{\csi}{2}
+\zeta(\xi)
+\left(\frac{-1+i \sqrt{3}}{2}\right)\left(\zeta(\xi-\xi_{i,0}^{\psi})+\zeta(\xi_{i,0}^{\psi})\right)
}\\ \displaystyle{ \phantom{1234567890123456789012345678901234567}
+\left(\frac{-1-i \sqrt{3}}{2}\right)\left(\zeta(\xi-\xi_{i,1}^{\psi})+\zeta(\xi_{i,1}^{\psi})\right),
}\\ \displaystyle{ 
\csi\not=0:\
\psi= - j \frac{9 \ex-4 j \ey}{24 \ex} \csi 
+\frac{j}{2} \left(\zeta(\xi-\xi_{j}^{\psi})+\zeta(\xi_{j}^{\psi})- \zeta(\xi))\right)
}\\ \displaystyle{ \phantom{12345678901234567890123}
+\frac{\sqrt{3}}{2} \left(\zeta(\xi-\xi_{j,0}^{\psi})+\zeta(\xi_{j,0}^{\psi})-\zeta(\xi-\xi_{j,1}^{\psi})-\zeta(\xi_{j,1}^{\psi})\right),
}\\ \displaystyle{
\csi\not=0:\
\frac{M'}{M} =  \frac{3 \ex+4 j \ey}{12 \ex} \csi 
+\zeta(\xi-\xi_{j}^{\psi})+\zeta(\xi_{j}^{\psi})+ \zeta(\xi)
}\\ \displaystyle{ \phantom{1234567890123456789}
-\left(\zeta(\xi-\xi_{j,0}^{\psi})+\zeta(\xi_{j,0}^{\psi})+\zeta(\xi-\xi_{j,1}^{\psi})+\zeta(\xi_{j,1}^{\psi}) )\right).
}
\end{array}
\right.
\label{eqsumzeta_G}
\end{eqnarray}
The choice $j=+i$ must be made for (\ref{eqdlogA}) to hold true,
while the choice $j=-i$ corresponds to the relation
\begin{eqnarray}
& & 
\frac{\D}{\D \xi} \log \left(\overline{A} e^{\displaystyle{-i \omega t + i \frac{c s_r}{2} \xi}}\right)
 =\frac{M'}{2 M} - i \psi.
\label{eqdlogAcc}
\end{eqnarray}

Before taking the quadrature of the above three expressions (\ref{eqsumzeta_G}),
let us recall the definition of the \textit{\'el\'ement simple} 
\cite[vol.~II, p.~506]{HalphenTraite} introduced by Hermite 
for integrating the Lam\'e equation,
\begin{eqnarray}
& &
\Elemsimp(\xi,q,k)=\frac{\sigma(\xi+q)}{\sigma(\xi) \sigma(q)} e^{(k-\zeta(q)) \xi},\
(k,q) \hbox{ constants}.
\label{eqdefElemsimp}
\end{eqnarray}
Its only singularity is a simple pole with residue unity at the origin.

Equations (\ref{eqsumzeta_G})${}_1$ and (\ref{eqsumzeta_G})${}_3$ then integrate as
\begin{eqnarray} 
& & {\hskip -18.0 truemm}
\forall \csi:\
A =K_0 
e^{\displaystyle{-i \omega t + i \frac{c \xi}{2 p}}} \
\Elemsimp(\xi,-\xi_{i,0}^{\psi},0)^{(-1+i \sqrt{3})/2} 
\Elemsimp(\xi,-\xi_{i,1}^{\psi},0)^{(-1-i \sqrt{3})/2},
\label{eqCGL5AprodHermite_G}
\\
& & {\hskip -18.0 truemm}
\csi\not=0:\
M =K_1 
e^{\displaystyle{\frac{3 \ex+4 j \ey}{12 \ex} \csi \xi}} \
\Elemsimp(\xi,-\xi_{j}^{\psi},0)
\Elemsimp(\xi,-\xi_{j,0}^{\psi},0)^{-1} 
\Elemsimp(\xi,-\xi_{j,1}^{\psi},0)^{-1}.
\label{eqCGL5MprodHermite_G}
\end{eqnarray}
The integration constants $K_0,K_1$ are determined
by requiring that, near the simple pole $\chi_1=\xi-\xi_{j,0}^{\psi}$ $\to 0$ of $M$,
the variables $A$ and $M$ admit the principal parts
$A \sim A_0   \chi_1^{(-1 + i \sqrt{3})/2}$,
$M \sim A_0^2 \chi_1^{-1}$, $A_0^8=3 / e_i^2$,
see (\ref{eqleadingA5}) and (\ref{eqCGL5_poles_M}). 

In order to check that the product of the complex amplitude $A$ (\ref{eqCGL5AprodHermite_G}) 
by its complex conjugate is equal to the decomposition (\ref{eqCGL5MprodHermite_G}),
one must take account of (\ref{eq_relations_among_poles}) and 
remember that the origin of $\xi$, not displayed in the above formulae,
depends on $j$ and is therefore different for $A$ and its complex conjugate.

The restriction $\csi\not=0$ in (\ref{eqCGL5MprodHermite_G})
is removed by taking into account the relation
\begin{eqnarray} 
& & 
 \zeta(\xi_{j,0}^{\psi} )+\zeta(\xi_{j,1}^{\psi} ) -\zeta(\xi_{j}^{\psi})
=\zeta(\xi_{j,0}^{\psi}  +      \xi_{j,1}^{\psi}   -      \xi_{j}^{\psi})
+\csi/4 +j \csi \ey/(3 \ex), 
\end{eqnarray}
and using the definition (\ref{eqEta1}), yielding
\begin{eqnarray} 
& & {\hskip -18.0 truemm} 
M=-K_1 
 e^{\displaystyle{ - \Eta_1 \xi} }  
\frac
{\sigma(\xi-\xi_{j  }^{\psi}) \sigma(\xi)}
{\sigma(\xi-\xi_{j,0}^{\psi}) \sigma(\xi-\xi_{j,1}^{\psi})} 
\frac
{\sigma(\xi_{j,0}^{\psi}) \sigma(\xi_{j,1}^{\psi})}
{\sigma(\xi_{j  }^{\psi})}.
\label{eqCGL5Mprodsigma_G}
\end{eqnarray}

In order to check the equality of the two decompositions of $M$ as
the sum (\ref{eqCGL5Msumzeta_g}) and the product  (\ref{eqCGL5Mprodsigma_G}),
it is sufficient to convert the elliptic function (\ref{eqCGL5Mprodsigma_G})
to a rational function of $\wp(\xi,G_2,G_3)$ and its derivative,
then to identify it to (\ref{eqCGL5Mwp}) \textit{modulo}
the Landen transformation (Appendix A).

Numerical simulations with periodic boundary conditions 
\cite[Fig.~4]{PSAK}
do display solutions $M$ having a real period
(similar features are observed in CGL3 \cite[Fig.~7]{Chate1994}),
these could well correspond to the present elliptic solution.

Remark.
The elliptic (hence singlevalued) nature of $\D \log (A e^{i \omega t})/\D \xi$
explains the so-called
``ad hoc Hirota method'' \cite{NB1984} 
in which $A$ is essentially assumed to be a product 
of powers of entire functions,
the powers being those of the singularity structure,
here $(-1\pm i \sqrt{3})/2$.
In order to recover our result (\ref{eqCGL5AprodHermite_G}),
two upgrades to this method are needed:
(i) to assume 
$A$ to be a product of powers
of Hermite's simple elements (\ref{eqdefElemsimp}),
not of Weierstrass $\sigma$ functions or Jacobi $\theta$ functions,
so as to ensure that the logarithmic derivative of $A$ is elliptic;
(ii) to allow arbitrary shifts $\xi_j$ in the arguments of the entire
functions, not only half periods like with the choice $\theta_j(\xi),j=0,1,2,3$
in the Jacobi notation. 

\appendix

\section*{Appendix A. Landen transformation}
\label{AppendixLanden}

We are indebted to the grateful indications of Yuri Brezhnev for this appendix.

The Landen or Gauss transformation consists in halving only one of the two periods,
it is naturally defined 
\cite{Kiepert1873} 
\cite[p.~384]{Kiepert1885}
in the notation of $\wp$ displaying the two periods $2 \omega,2 \omega'$,
\begin{eqnarray}
& & {\hskip -18.0 truemm}
 \wp(x        \vert  \omega,2\omega')
=\wp(x        \vert 2\omega,2\omega')
+\wp(x-\omega \vert 2\omega,2\omega')
-\wp(  \omega \vert 2\omega,2\omega').
\label{eqLanden_wp_omega}
\end{eqnarray}
In the other usual notation
\begin{eqnarray}
& & {\hskip -18.0 truemm}
\left\lbrace
\begin{array}{ll}
\displaystyle{
\wp  (x,G_2,G_3)= \wp(x        \vert  \omega,2 \omega'),\
\wp_1(x,g_2,g_3)= \wp(x        \vert 2\omega,2 \omega'),\
}\\ \displaystyle{
{\wp'  }^2=4 (\wp  ^3-g_2 \wp  -g_3)=4 (\wp  -e_1) (\wp  -e_2) (\wp  -e_3),\ 
}\\ \displaystyle{
{\wp_1'}^2=4 (\wp_1^3-G_2 \wp_1-G_3)=4 (\wp_1-E_1) (\wp_1-E_2) (\wp_1-E_3),\ 
}
\end{array}
\right.
\end{eqnarray}
the expression of $\wp(x \vert \omega,2\omega')$ 
as a rational function of $\wp(x \vert 2\omega,2\omega')$ is
\begin{eqnarray}
& & 
\wp(x,G_2,G_3)=\wp(x,g_2,g_3) -\frac{g_2-12 \ex^2}{4 (\wp(x,g_2,g_3)-\ex)},
\label{eqLanden_wp_g}
\end{eqnarray}
and similarly at the $\zeta$ and $\sigma$ levels \cite[Eqs.~(16b), (17b)]{Kiepert1885}
\begin{eqnarray}
{\hskip -18.0 truemm}
\zeta(x,G_2,G_3) &=& \zeta(x,g_2,g_3)+\zeta(x-\omega,g_2,g_3) - e_1 x + \zeta(\omega,g_2,g_3),
\label{eqLanden_zeta_g}
\\
\sigma(x,G_2,G_3) &=&
e^{\displaystyle{e_1 x^2/2 - \zeta(\omega,g_2,g_3) x}}\
\frac{\sigma(x,g_2,g_3)\sigma(x+\omega,g_2,g_3)}{\sigma(\omega,g_2,g_3)}.
\label{eqLanden_sigma_g}
\end{eqnarray}
Between 
$e_j,E_j$ (and $g_k,G_k$),
there exist two algebraic relations
\begin{eqnarray}
& & {\hskip -18.0 truemm}
\left\lbrace
\begin{array}{ll}
\displaystyle{
E_1=-2 e_1,\
(E_2-E_3)^2=36 e_1^2 - 4 (e_2-e_3)^2,
}\\ \displaystyle{
-32 g_2 g_3 + 22 g_3 G_2 + 11 g_2 G_3 - G_2 G_3=0,
}\\ \displaystyle{
196 g_2^3 + 49 g_2^2 G_2 -7260 g_3^2 + 660 g_3 G_3 - 15 G_3^2=0.
}
\end{array}
\right.
\label{eqLandenG}
\end{eqnarray}
The ratio $-2$ of the two zeros 
$-2 \ex$ of ${\wp_1'}^2$ in (\ref{eqCGL5psiwp})
and 
$\ex$ of ${\wp'}^2$ in (\ref{eqCGL5Mwp})
is the signature of such a Landen transformation.


\vfill \eject

\end{document}